\documentclass[10pt,leqno]{article}
\usepackage{amssymb,amsbsy,amsmath,amsfonts,amssymb,amscd, mathrsfs}

\usepackage[english]{babel}
\usepackage[T1]{fontenc}
\usepackage{indentfirst}

\makeatletter
\@addtoreset{equation}{section}
\makeatother

\newtheorem{statement}{}[section]
\newtheorem{theorem}[statement]{Theorem}
\newtheorem{lemma}[statement]{Lemma}

\newtheorem{proposition}[statement]{Proposition}

\newcommand\C{\mathbb C}

\newcommand\R{\mathbb R}
\newcommand\T{\mathbb T}
\newcommand\D{\mathbb D}

\newcommand\e{{\rm e}}

\newcommand\esp{\mathbb E}
\newcommand\eps{\varepsilon}
\newcommand\ind{{\rm 1\kern-.30em I}}
\newcommand\qed{\hfill $\square$}
\renewcommand \Re{{\mathfrak R}{\rm e}\,}
\renewcommand \Im{{\mathfrak I}{\rm m}\,}
\let\amphi=\phi
\let\phi=\varphi

\title{\bf Infinitesimal Carleson property for weighted measures induced by analytic self-maps of the unit disk}
\author{\it Daniel Li, Herv\'e Queff\'elec, Luis Rodr{\'\i}guez-Piazza}

\date{\footnotesize \today}

\begin{document}

\maketitle

\noindent{\bf Abstract.} \emph{We prove that, for every $\alpha > -1$, the pull-back measure $\phi ({\cal A}_\alpha )$ of the measure 
$d{\cal A}_\alpha (z) = (\alpha + 1) (1 - |z|^2)^\alpha \, d{\cal A} (z)$, where ${\cal A}$ is the normalized area measure on the unit disk $\D$, by 
every analytic self-map $\phi \colon \D \to \D$ is not only an $(\alpha + 2)$-Carleson measure, but that the measure of the 
Carleson windows of size $\eps h$ is controlled by $\eps^{\alpha + 2}$ times the measure of the corresponding window of size $h$. This means that 
the property of being an $(\alpha + 2)$-Carleson measure is true at all infinitesimal scales. We give an application by characterizing the 
compactness of composition operators on weighted Bergman-Orlicz spaces.}  
\medskip

\medskip


\medskip

{\scriptsize 

\noindent 
\it {\rm Daniel Li},  Univ Lille-Nord-de-France UArtois, 
Laboratoire de Math\'ematiques de Lens EA~2462, 
F\'ed\'eration CNRS Nord-Pas-de-Calais FR~2956, 
Facult\'e des Sciences Jean Perrin,
Rue Jean Souvraz, S.P.\kern 1mm 18, 
F-62\kern 1mm 300 LENS, FRANCE -- 
daniel.li@euler.univ-artois.fr
\smallskip

\noindent
{\rm Herv\'e Queff\'elec},   
Univ Lille-Nord-de-France USTL, 
Laboratoire Paul Painlev\'e U.M.R. CNRS 8524, 
F-59\kern 1mm 655 VILLENEUVE D'ASCQ Cedex, FRANCE --  
herve.queffelec@univ-lille1.fr
\smallskip

\noindent
{\rm Luis Rodr{\'\i}guez-Piazza},   
Universidad de Sevilla, 
Facultad de Matem\'aticas, Departamento de An\'alisis Matem\'atico,
Apartado de Correos 1160,
41\kern 1mm 080 SEVILLA, SPAIN --  
piazza@us.es\par}


\section{Introduction and notation} 

It is well-known that every analytic self-map $\phi \colon \D \to \D$ induces a bounded composition operator $f \mapsto C_\phi (f) = f \circ \phi$ 
from the Bergman space ${\mathfrak B}^2$ into itself. By Hastings's version of the Carleson inclusion theorem (\cite{Hastings}), that means that the 
pull-back measure ${\cal A}_\phi$ of the normalized area measure ${\cal A}$ by $\phi$ is a $2$-Carleson measure, that is, for some constant $C > 0$,  
\begin{displaymath}
{\cal A} \big( \{z\in \D\,;\ \phi (z) \in W (\xi, \eps) \} \big) \leq C\, \eps^2
\end{displaymath}
for every $\eps \in (0,1)$ and every $\xi \in \T$, where $W (\xi, \eps)$ is the Carleson window centered at $\xi$ and of size $\eps$. It was 
proved in \cite{Bergman}, Theorem~3.1, that one actually has an infinitesimal version of this property, namely, for some constant $C > 0$: 
\begin{equation}\label{contraction Bergman}
{\cal A} \big( \{z \in \D\,;\ \phi (z) \in W (\xi, \eps h) \} \big) 
\leq C  \, {\cal A} \big( \{z \in \D\,;\ \phi (z) \in W (\xi, h) \} \big) \,\eps^2 \,,
\end{equation}
for every $\eps \in (0, 1)$ and $h > 0$ small enough.
\medskip

Now, consider, for $\alpha > -1$, the \emph{weighted Bergman space} ${\mathfrak B}^2_\alpha$. By Littlewood's subordination principle, every 
analytic self-map $\phi$ of $\D$ induces a bounded composition operator $C_\phi$ from ${\mathfrak B}^2_\alpha$ into itself (see 
\cite{McCluer-Shapiro}, Proposition~3.4). By Stegenga's version of the Carleson theorem (\cite{Stegenga}, Theorem~1.2), that means that the 
pull-back measure of ${\cal A}_\alpha$ (see \eqref{mesure a poids} below) by $\phi$ is an $(\alpha + 2)$-Carleson measure. Our goal in this paper is 
to show the analog of \eqref{contraction Bergman} in the following form.
\par\medskip

\begin{theorem} \label{theo principal}
For each $\alpha > - 1$, there exists a constant $C_\alpha > 0$ such that, for every analytic self-map of the unit disk $\phi \colon \D \to \D$, 
every $\eps \in(0, 1)$ and every $h > 0$ small enough, one has, for every $\xi \in \T$:
\begin{equation} 
\begin{split}
{\cal A}_\alpha \big( \{z \in \D\,;\ \phi (z) \in W (\xi, & \eps h) \} \big)  \\
& \leq C_\alpha \,\eps^{\alpha + 2} \, {\cal A}_\alpha \big( \{z \in \D\,;\ \phi (z) \in W (\xi, h) \} \big)\,.
\end{split}
\end{equation}
\end{theorem}
\medskip

It should be stressed that the heart of the proof given in \cite{Bergman} in the case $\alpha = 0$ cannot be directly used for the other 
$\alpha > - 1$, and we have to change it, justifying the current paper. Moreover, the present proof is simpler than that of \cite{Bergman}. We also pointed out 
that the result holds in the limiting case $\alpha = -1$, corresponding to the Hardy space $H^2$ (\cite{CompOrli}, Theorem~4.19), but the proof is different, 
due to the fact that one uses the normalized Lebesgue measure on $\T$ and the boundary values of $\phi$ instead of measures on $\D$ and the function 
$\phi$ itself. 

\medskip

We end the paper by an application to the compactness of composition operators on weighted Bergman-Orlicz spaces.\par\smallskip

Another application of Theorem~\ref{theo principal} is given in \cite{approx}. 
\medskip

\noindent{\bf Notation.} In this paper, $\D = \{z\in \C\,;\ |z| < 1\}$ denotes the open unit disk of the complex plane $\C$, and 
$\T = \partial \D$ is the unit circle. The normalized area measure $\frac{dx\,dy}{\pi}$ is denoted by ${\cal A}$. \par
\medskip

For $\alpha > -1$, the \emph{weighted Bergman space} ${\mathfrak B}^2_\alpha$ is the space of all analytic functions 
$f (z) = \sum_{n=0}^\infty a_n z^n$ on $\D$ such that 
\begin{displaymath}
\| f \|_\alpha^2 := \int_\D | f (z)|^2 \,d{\cal A}_\alpha (z) < + \infty \,,
\end{displaymath} 
where ${\cal A}_\alpha$ is the weighted measure 
\begin{equation}\label{mesure a poids}
d{\cal A}_\alpha (z) = (\alpha + 1)(1 - |z|^2)^\alpha d{\cal A} (z)\,. 
\end{equation}

The \emph{Carleson window} centered at $\xi \in \T$ and of size $h$, $0 < h < 1$, is the set 
\begin{displaymath} 
W (\xi, h) = \{z \in \D\,;\ |z| \geq 1 - h \text{ and } |\arg (z \bar{\xi}) | \leq h \}\,.
\end{displaymath} 
A measure $\mu$ on $\D$ is called an \emph{$\alpha$-Carleson measure} ($\alpha \geq - 1$) if 
\begin{displaymath} 
\sup_{|\xi| = 1} \mu [ W (\xi, h)] = O_{h \to 0} \,(h^\alpha).
\end{displaymath} 
Actually, instead of the Carleson window $W (\xi, h)$, we shall merely use the sets 
\begin{displaymath} 
S (\xi, h) = \{z \in \D\,;\ |z - \xi | \leq h \}\,,
\end{displaymath} 
which have essentially the same size, so $\mu$ is an \emph{$\alpha$-Carleson measure} if and only if 
$\sup_{|\xi| = 1} \mu [ S (\xi, h)] = O_{h \to 0} \,(h^\alpha)$.\par
\medskip

We denote by $\Pi^+$ the right-half plane 
\begin{equation} \label{half plane}
\Pi^+ = \{z \in \C\,;\ \Re z > 0 \} \,.
\end{equation} 
To avoid any misunderstanding, we denote by $A$ the area measure on $\Pi^+$, and \emph{not} this measure divided by $\pi$. \par\smallskip

Let $T \colon \D \to \Pi^+$ be the conformal map defined by: 
\begin{equation} \label{transfo conforme}
T (z) = \frac{1 - z}{1 + z} \,;
\end{equation} 
we denote by $\tau_\alpha = T ({\cal A}_\alpha)$ the pull-back measure defined by:
\begin{equation} 
\tau_\alpha (B) = {\cal A_\alpha} [T^{- 1} (B)]
\end{equation} 
for every Borel set $B$ of $\Pi^+$. This is a probability measure on $\Pi^+$.  \par\smallskip

We also need another measure $\mu_\alpha$ on $\Pi^+$, defined by:
\begin{equation} 
d\mu_\alpha = x^\alpha \,dx dy \,.
\end{equation} 

Given two measures $\mu$ and $\nu$, we shall write $\mu \sim \nu$ when the Radon-Nikod\'ym derivative $\frac{d \mu}{d \nu}$ is bounded from above 
and from below. \par
\medskip

The \emph{pseudo-hyperbolic distance} $\rho'$ on $\D$ is given by
\begin{equation} 
\rho' (z, w) = \bigg| \frac{z - w}{1 - \bar{z} w} \bigg| \raise 1pt \hbox{,} \qquad z, w \in \D\,.
\end{equation} 
For every $z \in \D$ and $r \in (0, 1)$, 
\begin{displaymath} 
\Delta' (z, r) = \{w \in \D \, ; \ \rho' (w, z) < r \} 
\end{displaymath} 
is called the \emph{pseudo-hyperbolic disk} with center $z$ and radius $r$. It is (see \cite{BM}, \cite{Duren-S}, or \cite{Zhu}, for example) the image of the 
Euclidean  disk $D (0, r)$ by the automorphism 
\begin{displaymath} 
\phi_z (\zeta) = \frac{z - \zeta}{1 - \bar{z} \zeta} \, \cdot
\end{displaymath} 
\par

The pseudo-hyperbolic distance $\rho$ on $\Pi^+$ is deduced by transferring the pseudo-hyperbolic distance $\rho'$ on $\D$ with the conformal map $T$: 
\begin{equation} 
\rho (a, b) = \rho' ( T^{-1} a, T^{-1} b) = \bigg| \frac{a - b}{\bar{a} + b}\bigg| \, \raise 1pt \hbox{,}
\end{equation} 
and, for every $w \in \Pi^+$ and $r \in (0, 1)$, 
\begin{displaymath} 
\Delta (w, r) = \{z \in \Pi^+ \, ; \ \rho (z, w) < r \} 
\end{displaymath} 
is the \emph{pseudo-hyperbolic disk} of $\Pi^+$ with center $w$ and radius $r$. 
\medskip

Finally, we shall use the following notation:
\begin{equation} \label{Omega}
\Omega \, = ( 0 , 2 ) \times ( - 1 , 1 ) \,.
\end{equation} 
\medskip

\noindent{\bf Acknowledgement.} Part of this work was made during a visit of the second named author at the Departamento de An\'alisis Matem\'atico of 
the Universidad de Sevilla in April 2011; it is a pleasure to thank all people of this department for their warm hospitality. The third-named author is partially 
supported by a Spanish research project MTM 2009-08934. 
\goodbreak


\section{Transfer to the right half plane}

As in \cite{Bergman}, we only have to give the proof for $\xi = 1$ and, by considering $g = h/ (1 - \phi)$, we are boiled down to prove:

\begin{theorem} \label{reduction}
Let $\alpha > - 1$. There exist constants $K_0 > 0$, $c_0 > 0$ and $\lambda_0 > 1$ such that every analytic function $g \colon \D \to \Pi^+$ 
with $|g (0)| \leq c_0$ satisfies, for every $\lambda \geq \lambda_0$:
\begin{displaymath} 
{\cal A}_\alpha (\{ |g| > \lambda \}) \leq \frac{K_0}{\lambda^{\alpha + 2}}\, {\cal A}_\alpha (\{ |g| > 1\})\,.
\end{displaymath} 
\end{theorem} 

As said in the Introduction, this result is an infinitesimal version of the fact that the pull-back measure ${\cal A}_{\alpha, \phi}$ of 
${\cal A}_\alpha$ by any analytic self-map $\phi$ of $\D$ is an $(\alpha + 2)$-Carleson measure. In fact, one has the following result.

\begin{proposition} \label{starting lemma}
There is some constant $C = C_\alpha > 0$ such that
\begin{equation} \label{eq 1 starting lemma}
{\cal A}_\alpha (\{ |g| > \lambda\}) \leq \frac{C}{\lambda^{\alpha + 2}}\, |g (0)|^{\alpha + 2}
\end{equation} 
for every analytic function $g \colon \D \to \Pi^+$ and every $\lambda > 0$. 
\end{proposition}

The goal is hence to replace in the right-hand side the quantity $|g (0)|^{\alpha + 2}$ by ${\cal A}_\alpha (\{ |g| > 1\})$. \par
\bigskip

\noindent{\bf Proof of Proposition~\ref{starting lemma}.} We may assume that $| g (0)| = 1$. Hence we may assume that $\lambda > 2$, taking 
$C \geq 2^{\alpha + 2}$, because ${\cal A}_\alpha (\{ |g| > \lambda\}) \leq 1$.\par\smallskip
Set $\phi (z) = [g (z) - g (0)]/[g (z) + \overline{g (0)}]$. Then $|g (z)| > \lambda$ implies that 
\begin{displaymath} 
|\phi (z) - 1| = 2 \, |\Re g (0)|/ | g (z) + \overline{g (0)}| \leq 2/(\lambda - 1) \leq 4/ \lambda\,.
\end{displaymath} 
But $\phi$ maps $\D$ into itself, so the measure ${\cal A}_{\alpha, \phi}$ is an $(\alpha + 2)$-Carleson measure and (see the proof of 
\cite{Stegenga}, Theorem~1.2) 
\begin{displaymath} 
{\cal A}_\alpha (\{ |g| > \lambda\}) \leq {\cal A}_{\alpha, \phi} [ S (1, 4/\lambda)] \leq C'_\alpha \,\| C_\phi \|^2 / (\lambda/4)^{\alpha + 2} \,,
\end{displaymath} 
where $\| C_\phi \|$ is the norm of the composition operator $C_\phi \colon {\mathfrak B}_\alpha^2 \to {\mathfrak B}_\alpha^2$. But $\phi (0) = 0$ 
and hence $\| C_\phi \| = 1$, by using Littlewood's subordination principle and integrating. 
\qed \par
\bigskip

For technical reasons, that we shall explain after Lemma~\ref{touch}, we need to work with functions defined on $\Pi^+$. 
Proposition~\ref{starting lemma} becomes:
\begin{proposition} \label{transfer}
There exists a constant $C = C_\alpha > 0$ such that, for every analytic function $f \colon \Pi^+ \to \Pi^+$, one has: 
\begin{equation} \label{global} 
\tau_{\alpha} (\{ \vert f \vert > \lambda \} ) \leq \frac{C}{\lambda^{\alpha + 2}} \vert f (1) \vert^{\alpha+2}. 
\end{equation}
\end{proposition}

\noindent{\bf Proof.} Set $E_{f} (\lambda) = \{\vert f \vert > \lambda\}$ and define similarly $E_{g} (\lambda) = \{\vert g \vert > \lambda\}$ where 
$g = f \circ T \colon \D \to \Pi^+$. We have $g (0) = f (1)$ as well as the simple but useful equation:
 \begin{equation} \label{useful}
T^{-1}[E_{f} (\lambda)] = E_{g} (\lambda).
\end{equation} 
So that, by Proposition~\ref{starting lemma}:
\begin{align*}
\tau_{\alpha}[E_{f} (\lambda)] 
& = {\cal A}_{\alpha} [T^{-1} (E_{f} (\lambda) ] = {\cal A}_{\alpha} [E_{g} (\lambda) ] \\ 
& \leq \frac{C}{\lambda^{\alpha + 2}} \vert g (0) \vert^{\alpha+2} 
= \frac{C}{\lambda^{\alpha+2}} \vert f (1) \vert^{\alpha+2} \,,
\end{align*}
and Proposition~\ref{transfer} is proved. 
\qed
\bigskip

Now, to prove Theorem~\ref{reduction}, it suffices to prove that, when one localizes $f$ on $\Omega$, one may replace the quantity $|f (1)|$ in the 
right-hand side of \eqref{global} by $\tau_\alpha ( \{|f| > 1\} \cap \Omega)$. This is what is claimed in the next result. 

\begin{theorem} \label{theo clef}
There exist constants $K = K_\alpha > 0$, $c_1 > 0$ and $\lambda_1 > 1$ such that every analytic function $f \colon \Pi^+ \to \Pi^+$ such that 
$| f ( 1)| \leq c_1$ satisfies, for every $\lambda \geq \lambda_1$: 
\begin{displaymath} 
\tau_{\alpha} (\{ \vert f \vert > \lambda\} \cap\Omega ) 
\leq \frac{K}{\lambda^{\alpha + 2}} \, \tau_{\alpha} (\{ \vert f \vert > 1\} \cap \Omega).
\end{displaymath} 
\end{theorem} 
\bigskip

We shall prove Theorem~\ref{theo clef} in the next section, but before, let us see why it gives Theorem~\ref{reduction} and hence our main result,  
Theorem~\ref{theo principal}.
\bigskip
\goodbreak

\noindent{\bf Proof of Theorem~\ref{reduction}.} Let $E \colon \Pi^+ \to \D$ be the exponential map defined by 
\begin{equation} 
E (z) = \e^{- \pi z} \,,
\end{equation} 
which (up to a radius) maps bijectively $\Omega$ onto the annulus
\begin{equation} \label{annulus}
U = \{z \in \D \, ; \ \vert z \vert > \e^{- 2\pi}\} . 
\end{equation} 

For every $g \colon \D \to \Pi^+$ with $|g (0)| \leq (1 - \beta)/ (1 + \beta)$ and $0 < \beta < 1$, one has, by Schwarz's lemma (see \cite{Bergman}, eq.~(3.9)):
\begin{displaymath} 
|g (z)| > 1 \qquad \Longrightarrow \qquad |z| > \beta \,.
\end{displaymath} 
Therefore we only have to work on the annulus $U$, taking $c_0 \leq \tanh \pi$ in Theorem~\ref{reduction}. \par
\medskip

Let $L = E^{-1}$ be the inverse map of the restriction of $E$ to $\Omega$, and 
\begin{equation} 
\sigma_\alpha = L ({\cal A}_\alpha) 
\end{equation} 
be the pull-back measure of ${\cal A}_\alpha$ by $L$. This measure is carried by $\Omega$ and we have: 
\begin{lemma} \label{lastminute} 
On $\Omega$, one has: $\sigma_\alpha\sim \mu_\alpha\sim \tau_\alpha$.
\end{lemma}
\medskip

Taking this lemma for granted for a while, let us finish the proof of Theorem~\ref{reduction} (the measure $\mu_\alpha$ does not come into play here). Let 
$g \colon \D \to \Pi^+$ be an analytic function and $f = g \circ E \colon \Pi^+ \to \Pi^+$ (so that $g = f \circ L$ on $E (\Omega)$). We have 
$| f (1) | \leq c_1$ if $| g (0) | \leq c_0$, with $c_0 > 0$ small enough. In fact, the analytic function $h = T \circ g$ maps $\D$ into itself and hence, by the 
Schwarz-Pick inequality, $h$ is a contraction for the pseudo-hyperbolic distance on $\D$ (see \cite{BM}, eq.~(3.3), page~18, for example);  
hence $\rho' [h (\e^{-\pi}), h (0)] \leq \rho' (\e^{-\pi}, 0) = \e^{-\pi}$, that is 
$\Big| \frac{g (\e^{-\pi}) - g (0)}{g (\e^{-\pi}) + \overline{g (0)}} \Big| \leq \e^{-\pi}$. It follows that 
$|g (\e^{-\pi})| - | g (0)| \leq \e^{-\pi} \big[ | g (\e^{-\pi}) | + | g (0)| \big]$, {\it i.e.} $|g ( \e^{-\pi}) | \leq \frac{1}{\tanh \pi} \, |g (0)|$. Therefore 
$|f (1)| = | g (\e^{-\pi})| \leq c_1$ if $| g (0)| \leq c_0$ with $c_0 \leq c_1 \tanh \pi$.

Set:
\begin{displaymath} 
\qquad E_{g} (\lambda) = \{\vert g \vert>\lambda\} \cap U \qquad  \text{and} \qquad E_{f} (\lambda) = \{\vert f \vert>\lambda\}\cap \Omega .
\end{displaymath} 
Observe that, as in \eqref{useful}, 
\begin{displaymath} 
\qquad  L^{-1} [E_{f} (\lambda)] = E_{g} (\lambda) \qquad \text{and}\qquad  E^{-1}[E_{g} (1)] = E_{f} (1) .  
\end{displaymath} 
Hence, in view of Theorem~\ref{theo clef} and Lemma~\ref{lastminute}:
\begin{align*}
{\cal A} _{\alpha} [E_{g} (\lambda)] 
& = {\cal A}_{\alpha} \big(L^{-1} [E_{f}(\lambda)] \big) 
=\sigma_{\alpha} [E_{f} (\lambda)] \\
& \leq \frac{K'_\alpha}{\lambda^{\alpha+2}} \, \sigma_{\alpha} [E_{f} (1) ] 
= \frac{K'_\alpha}{\lambda^{\alpha+2}} \, \sigma_{\alpha}\big(E^{-1} [E_{g} (1)] \big) \\
& = \frac{K'_\alpha}{\lambda^{\alpha+2}} \, (E\sigma_{\alpha}) [E_{g}(1)] 
= \frac{K'_\alpha}{\lambda^{\alpha+2}} \, {\cal A}_{\alpha} [E_{g}(1)] \,,
\end{align*}
which is exactly what we wanted to prove. 
\qed
\bigskip\goodbreak

\noindent{\bf Proof of Lemma~\ref{lastminute}.} Let us compute  $\sigma_\alpha$ with the change of variable $w = E^{-1} (z)$. One has $z = E (w)$ and 
\begin{displaymath} 
d{\cal A} (z) = \vert E' (w) \vert^2 \, \frac{dA (w)}{\pi} = \frac{1}{\pi}\e^{- 2\pi \Re w} \, dA (w).
\end{displaymath} 
We get: 
\begin{align*}
\int_{\Omega} h (w) \,d\sigma_{\alpha} (w) 
& =\int_{U} h (L z) \, d{\cal A}_{\alpha} (z) 
= (\alpha + 1) \int_{U} h (E^{-1} z) (1 - \vert z \vert^2)^\alpha d{\cal A} (z) \\
& = \frac{\alpha + 1}{\pi} \int_{\Omega} h (w) \, \e^{-2 \pi \Re w} (1 - \e^{-2 \pi \Re w})^\alpha \, d A(w) , 
\end{align*}
so that 
\begin{equation} \label{equiv sigma}
d \sigma_\alpha (w) = \frac{\alpha + 1}{\pi} \, \e^{-2 \pi \Re w} (1 - \e^{-2 \pi \Re w})^\alpha \, \ind_\Omega (w) \,d A (w) .
\end{equation} 
Thus, on $\Omega$, we have  $\sigma_\alpha \sim \mu_\alpha$.  
Indeed, the factor $\e^{-2 \Re w}$ is bounded from below and from above, and $(1 - \e^{- 2 \Re w})^\alpha\sim (\Re w)^\alpha$ as $\Re w$ goes to $0$. 
This proves the first equivalence of Lemma~\ref{lastminute}.  \par
\smallskip

To prove the second equivalence, we use the change of variable formula $z = T w$ in 
\begin{displaymath} 
\int_\Omega h (u)\, d\tau_\alpha (u) = \int_U h (T z) \, d{\cal A}_\alpha (z) ; 
\end{displaymath} 
it gives $d\tau_\alpha (w) = \vert T ' (w) \vert^2 (1 - \vert T (w) \vert^2)^\alpha (\alpha + 1) \, d A (w) / \pi$, i.e.:
\begin{equation} \label{image} 
d\tau_\alpha (w) = \frac{4^{\alpha + 1} (\alpha + 1)}{\pi} \, \frac{(\Re w)^\alpha}{\vert 1 + w \vert^{2 (\alpha + 2)}}\, \ind_\Omega (w) \, d A (w) ,
\end{equation}
showing that $\mu_\alpha \sim \tau_\alpha$ on $\Omega$. 
\qed


\section{Proof of Theorem~\ref{theo clef}} 

Let us split, up to a set of measure $0$, the square  $\Omega$  into dyadic sub-squares 
\begin{equation} 
Q_l = \bigg( \frac{2j}{2^n} \, \raise 1pt \hbox{,} \, \frac{2 (j+1)}{2^n} \bigg) 
\times \bigg( \frac{2k}{2^n} - 1 \, \raise 1pt \hbox{,} \, \frac{2 (k+1)}{2^n} - 1 \bigg)
\end{equation} 
 of center 
\begin{equation} 
c_l =  \frac{2j + 1}{2^n} + i\, \bigg(\frac{2k + 1}{2^n} - 1 \bigg),  
\end{equation} 
with $n \geq 0$ , $0 \leq j, k \leq 2^{n} - 1$ and where $l = (n, j, k)$.\par\smallskip

Note that $\Omega = Q_{(0, 0, 0)}$. We are going to use the special form of the measure $\tau_\alpha$, taken in \eqref{image}, to get a localized version of 
Proposition~\ref{transfer} as follows.
\goodbreak

 \begin{proposition} \label{auxiliary} 
There is a constant $C_\alpha > 0$ such that, for any analytic function $f \colon \Pi^+ \to \Pi^+$ and any dyadic sub-square $Q_l$ of $\Omega$, one has, 
for any $\lambda > 0$:
\begin{equation} \label{localize} 
\tau_{\alpha} ( \{\vert f \vert > \lambda\} \cap Q_l ) 
\leq \frac{C_\alpha}{\lambda^{\alpha+2}} \, \tau_{\alpha} (Q_l) \, \vert f (c_l) \vert^{\alpha + 2} . 
\end{equation}
\end{proposition}

\noindent{\bf Proof.} Using Lemma~\ref{lastminute}, we may replace the measure  $\tau_\alpha$ by $d \mu_\alpha = x^\alpha \,dx dy$. This measure is 
no longer a probability measure, but it has the advantage of being invariant under vertical translations, and, especially, to react to a dilation of positive ratio 
$\lambda$ by multiplying the result by the factor $\lambda^{\alpha + 2}$.  \par \smallskip

We first need a simple lemma. 

\begin{lemma} \label{harnack}  
For every $0 \leq s < 1$, there exists a constant $M_s > 0$ such that, for any analytic function $f \colon \Pi^+ \to \Pi^+$ and any pseudo-hyperbolic disk 
$\Delta (c, s)$ in $\Pi^+$, we have, for every $z \in \Delta (c, s)$: 
\begin{equation} \label {schwarz} 
1/ M_s \leq \vert f (z) \vert / \vert f (c) \vert \leq M_s.
\end{equation}
\end{lemma}

\noindent{\bf Proof.} By the classical Schwarz-Pick inequality, any analytic map $f \colon \Pi^+\to \Pi^+$ contracts the  pseudo-hyperbolic distance $\rho$ of 
$\Pi^+$ (see \cite{BM}, Section~6), so that if $z \in \Delta (c, s)$, one has:
\begin{displaymath} 
\vert u \vert := \bigg\vert \frac{f (z) - f (c)}{f (z) + \overline{f (c)}} \bigg\vert  \leq \bigg\vert \frac{z - c}{z + \overline{c}} \bigg\vert \leq s. 
\end{displaymath} 
Inverting that relation, we get $f (z) = \frac{u\overline{f (c)}+f(c)}{1 - u}$, whence 
\begin{displaymath} 
| f (z) | \leq | f (c) | \, \frac{1 + | u |}{1 - | u  |} \leq | f (c) | \, \frac{1 + s}{1 - s} 
\end{displaymath} 
and, similarly, $| f (z) | \geq | f (c) | \, \frac{1 - s}{1 + s}\,$. The lemma follows, with $M_s =  \frac{1 + s}{1 - s}\, \cdot$ 
\qed
\par\medskip

Let us now continue the proof of  Proposition~\ref{auxiliary}.\par
 
\begin{lemma} \label{notouch}  
Inequality \eqref{localize} holds when the square $Q_l$, of the $n$-th generation, does not touch the boundary of  $\Pi^+$, namely when $l = (n, j, k)$ with 
$j \geq 1$. More precisely, we have $Q_l \subseteq \Delta (c_l, s)$ where $s < 1$ is a numerical constant.
\end{lemma}

\noindent{\bf Proof.} Recall that $c_l$ is the center  of $Q_l$. We claim that we can find some numerical $s < 1$ such that $Q_l \subset \Delta (c_l, s)$. 
To show that claim, let $l = (n, j, k)$ and  $z, w \in Q_l$. We  have: 
\begin{displaymath} 
1 - \rho (z, w)^2 = 1 - \bigg\vert \frac{z - w}{z + \bar{w}} \bigg\vert^2 
= 4\, \frac{\Re z \, \Re w}{\vert z + \bar{w}\vert^2} \,\cdot
\end{displaymath} 
But one has $2j /2^{n} \leq \Re z, \Re w \leq 2 (j + 1)/ 2^n$ whereas $\vert \Im (z + \bar{w}) \vert \leq  2^{- n + 1}$; hence 
$\Re z \, \Re w \geq 4 j^2 4^{- n}$ and 
$\vert z + \bar{w}\vert^2 = (\Re z + \Re w)^2 + [\Im (z + \bar{w})]^2 \leq 16 (j + 1)^2 4^{- n} + 4. 4^{- n} \leq 80 j^2 4^{- n}$, because $j \geq 1$. 
Therefore
\begin{displaymath} 
1 - \rho (z, w)^2 \geq 4\, \frac{4 j^2 4^{- n}}{80 j^2 4^{- n}} = \frac{1}{5} \,\raise 1pt \hbox{,}
\end{displaymath} 
so that $\rho (z, w) \leq s =\sqrt{4/5}$. In particular, we have  $Q_l \subseteq \Delta (c_l, s)$. 
\par\medskip\goodbreak

Now, to prove \eqref{localize}, we may assume, by homogeneity (replace $f$ by $f/ |f (c_l) |$ and $\lambda$ by $\lambda / |f (c_l)|$), that $| f (c_l) | = 1$. 
We then have, by Lemma~\ref{harnack}, $\vert f (z) \vert \leq M_s \vert f (c_l) \vert = M_s$ for every $z \in Q_l$. Hence \eqref{localize} trivially holds when 
$\lambda > M_s$, since then the set in the left-hand side is empty. So we assume $\lambda \leq M_s$. In that case, setting $C_\alpha = M_s^{\alpha + 2}$, 
we have : 
\begin{displaymath} 
\tau_{\alpha} (\{\vert f \vert >\lambda\} \cap Q_l ) \leq \tau_{\alpha} (Q_l) \leq \frac{C_\alpha}{\lambda^{\alpha + 2}}\, \tau_{\alpha} (Q_l). 
\end{displaymath} 
This is the desired inequality, since we have supposed that $| f (c_l) | = 1$. 
\qed

\begin{lemma} \label{touch} 
Inequality \eqref{localize} holds when the square $Q_l$, of the $n$-th generation, touches the boundary of  $\Pi^+$, namely when $l = (n, j, k)$ with $j = 0$.
\end{lemma}

\noindent{\bf Proof.} This case uses the specific properties of the measure $\mu_\alpha$. In view of Lemma~\ref{lastminute}, we have to prove that:
\begin{equation} \label{withmu} 
\mu_{\alpha}(\{\vert f \vert > \lambda\} \cap Q_l )\leq \frac{C_\alpha}{\lambda^{\alpha+2}}\, \mu_{\alpha} (Q_l) \, \vert f (c_l) \vert^{\alpha+2} , 
\end{equation} 
when the square $Q_l \subseteq \Omega$ is supported by the imaginary axis. We may again assume that $| f (c_l) | = 1$, and we proceed in three steps. \par
\smallskip

1) First, \eqref{withmu} holds if  $Q_l = Q_{(0, 0, 0)} = \Omega$: this is just what we have proved in Proposition~\ref{transfer} with \eqref{global}. 
\par\smallskip

2) For $h > 0$, \eqref{withmu} holds when $Q_l = h \Omega = ( 0, 2h ) \times (- h, h )$ is a square meeting the imaginary axis in an interval 
$(- h, h )$ centered at $0$. Indeed, setting $E_{f} (\lambda) = \{\vert f \vert > \lambda\}$ as well as $f_{h} (z) = f (hz)$, we easily check that 
\begin{equation} \label{sameset} 
E_{f} (\lambda) \cap h\,\Omega = h\, [E_{f_h} (\lambda) \cap \Omega ]. 
\end{equation} 
For example, if $v \in E_{f_h} (\lambda) \cap \Omega$, one has $\vert f (h v) \vert > \lambda$ and hence $w = hv \in E_{f} (\lambda) \cap h\Omega$, 
giving one inclusion in \eqref{sameset}; the other is proved similarly. Using the already mentioned $(\alpha + 2)$-homogeneity of the measure $\mu_\alpha$, we 
obtain, using \eqref{global} for $f_h$:
\begin{align*}
\mu_\alpha [E_{f}(\lambda) \cap h \Omega] 
& = \mu_\alpha [h \big(E_{f_h} (\lambda) \cap \Omega \big) ] 
= h^{\alpha + 2} \mu_\alpha [E_{f_h} (\lambda) \cap \Omega] \\ 
& \leq h^{\alpha + 2} \frac{C_\alpha}{\lambda^{\alpha + 2}} \, \vert f_{h}(1) \vert^{\alpha+2}
= \mu_{\alpha}(Q_l) \, \frac{C'_\alpha}{\lambda^{\alpha + 2}} \, \vert f (c_l) \vert^{\alpha+2} ,
\end{align*}
with $C'_\alpha = 4^{- (\alpha + 2)} (\alpha + 1) C_\alpha$, since the center $c_l$ of  $Q_l = h \Omega$ is $c_l = h$. \par\smallskip

3) Finally, \eqref{withmu} holds if $Q_l$ is any square supported by the imaginary axis. Indeed, this $Q_l$ is a vertical translate of the second case,  and the 
measure $\mu_\alpha$ is invariant under vertical translations, which exchange centers.\par
\smallskip

This ends the proof of the crucial Lemma~\ref{touch} and thereby that of Proposition~\ref{auxiliary}. 
\qed
\par\medskip\goodbreak

\noindent{\bf Remark.} We see here why it is better to work with functions $f \colon \Pi^+ \to \Pi^+$ instead of functions $g \colon \D \to \Pi^+$; if the 
invariance of $\mu_\alpha$ under vertical translations corresponds to the rotation invariance of ${\cal A}_\alpha$, the homogeneity of $\mu_\alpha$, used 
in part 2) of the proof, corresponds to an invariance by the automorphisms $\phi_a$ of $\D$, with real $a \in \D$, which is not shared by ${\cal A}_\alpha$, 
and writing a measure equivalent to ${\cal A}_\alpha$ having these properties is not so simple.
\par\bigskip

In order to exploit this proposition, we need the following precisions.
\begin{lemma} \label{precision} 
There exist constants $c > 0$ and  $\delta_0 > 0$, depending only on $\alpha$, such that for every $l$, there exists $R_l \subseteq Q_l$ with 
$\tau_\alpha (R_l) \geq c\, \tau_\alpha (Q_l)$ and, for every analytic map $f \colon \Pi^+ \to \Pi^+$,
\begin{equation} \label{proportion} 
\qquad  \quad \vert f (z) \vert > \delta_0 \, \vert f (c_l) \vert  \quad \text{for every } z \in R_l.
\end{equation}
\end{lemma}

\noindent{\bf Proof.} By Lemma~\ref{lastminute}, it suffices to prove this lemma with $\mu_\alpha$ instead of $\tau_\alpha$. Let us consider two cases.
\par\smallskip

1) If $l = (n, j, k)$ with $j \geq 1$, we can simply take $R_l = Q_l$, in view of Lemma~\ref{harnack} and Lemma~\ref{notouch}. \par\smallskip

2) If $l = (n, j, k)$ with $j = 0$, we may assume that $Q_l = \Omega = ( 0, 2 ) \times (- 1, 1 )$, since either vertical translations or dilations of positive ratio 
are isometries for the pseudo-hyperbolic distance on $\Pi^+$ and, on the other hand, multiply the $\mu_\alpha$-measure by $1$ or $h^{\alpha + 2}$ 
respectively. It follows that $c_l = 1$. We are going to check that $\Delta (1, 1/4) \subseteq \Omega = Q_l $, so that we can take $R_l = \Delta (c_l, 1/4) $. 
Indeed, set $t = 1/4$; if $\vert u \vert := \big\vert \frac{z - 1}{z + 1} \big\vert \leq  t$, we have  $z = \frac{1 + u}{1 - u}$ and 
\begin{align*} 
& 0 < \Re z  = \frac{1 - \vert u \vert^2}{\vert 1 - u \vert^2} \leq \frac{1 + t}{1 - t} < 2\,; \\
& \vert \Im z \vert  = \frac{2\, \vert \Im u \vert}{\vert 1 - u \vert^2} \leq \frac{2t}{(1 - t)^2} = \frac{8}{9} < 1.
\end{align*} 
Moreover, in view of Lemma~\ref{harnack}, \eqref{proportion} holds with $\delta_0 = M_{t}^{-1} = 3/5$. \par \smallskip   

Finally, the claim on the measures holds with $c = \mu_\alpha [\Delta (1, 1/4)] / \mu_\alpha (\Omega)$. 
\qed
\bigskip

Now, we want to control mean values of $f$ on some of the $Q_l$'s. In order to get that, we have to do a Calder\'on-Zygmund decomposition. \par\smallskip

To that end, we need to know that the mean of $|f|$ on $\Omega$ is small, namely less than $1$, if $| f (1)|$ is small enough. This is the aim of the next 
proposition. 
\begin{proposition} \label{no beans}
There exists a constant $C > 0$ such that, for every analytic function $f \colon \Pi^+ \to \Pi^+$, one has:
\begin{equation} \label{moyenne Omega}
|f (1) | \leq \iint_\Omega | f( x + iy)| \, \frac{dx dy}{\pi} \leq C\, | f (1)| \,.
\end{equation} 
Moreover, if $c$ is the center of an open square $Q$ contained in $\Pi^+$, then:
\begin{equation} \label{moyenne Q}
\frac{\pi}{4} \, |f (c)| \leq \frac{1}{A (Q)} \int_Q | f (z)| \, dA (z) \leq C\, \frac{ \pi}{4}\, |f (c)| \,.
\end{equation} 
\end{proposition}
\noindent{\bf Proof.} Let us see first that \eqref{moyenne Q} follows from \eqref{moyenne Omega}. Let $c = a + ib$ ($a > 0$ and $b \in \R$) be the center of 
the square $Q = ( a - h, a + h ) \times ( b - h, b + h )$, with $0 < h \leq a$. Consider the function $f_1$ defined by:
\begin{displaymath} 
\qquad \qquad f_1 (z) = f [\amphi (z) ] \,, \quad \text{where} \quad \amphi (z) = h z - h + a +ib \,.
\end{displaymath} 
Observe that $\amphi \colon \Pi^+ \to \Pi^+$ is an affine transformation sending $1$ onto $c$ and that $\amphi (\Omega) = Q$. Applying 
\eqref{moyenne Omega} to $f_1$ gives:
\begin{displaymath} 
\frac{\pi}{4} \, |f_1 (1)| \leq \frac{1}{A (\Omega)} \int_\Omega | f_1 (z)| \, dA (z) \leq C\, \frac{ \pi}{4}\, |f_1 (1)| \,.
\end{displaymath} 
This yields \eqref{moyenne Q} using an obvious change of variable and $f_1 (1) =  f (c)$. \par \smallskip

The left-hand side inequality in \eqref{moyenne Omega} is due to subharmonicity: consider the open disk $D$ of center $1$ and radius $1$; then 
$D \subseteq \Omega$ and, $|f|$ being subharmonic, we have:
\begin{displaymath} 
|f (1) | \leq \frac{1}{\pi} \iint_D | f( x + iy)| \, dx dy \leq \frac{1}{\pi}\iint_\Omega | f( x + iy)| \, dx dy \,.
\end{displaymath} 

We now prove the right-hand side inequality. Using Lemma~\ref{transfer} and the fact that $\mu_0 \sim \tau_0$ on $\Omega$ (note that $\mu_0 $ is just 
the area measure $A$ on $\Pi^+$), we have the existence of a constant $\kappa > 0$ such that, for all $\lambda > 0$:
\begin{equation} \label{kappa}
\mu_0 (\{ |f| > \lambda \} \cap \Omega) \leq \frac{\kappa}{\lambda^2} \, | f (1)|^2 \,.
\end{equation} 
From this estimate \eqref{kappa}, we can control the integral of $|f|$ over $\Omega$ (recall that $\mu_0 (\Omega) = 4$):
\begin{align*} 
\int_\Omega |f| \, d\mu_0 
& = \int_0^{+\infty} \!\!\! \mu_0 (\{ |f| > \lambda \} \cap \Omega) \, d\lambda \\
& \leq 4\, | f( 1)| + \int_{|f (1)|}^{+\infty} \frac{\kappa \, | f (1)|^2}{\lambda^2}\, d\lambda = (4 + \kappa)\, |f (1)| .
\end{align*} 
The proposition follows. 
\qed
\medskip

\noindent{\bf Remark.} We do not know if the constant $\pi/4$ in the left-hand side of \eqref{moyenne Q} can be replaced
by a better constant; however, it is not possible to replace this factor $\pi/4$ by $1$. Let us see an example.\par

Let us define $f (z) = \exp \bigl( (Tz)^4) \bigr)$ where $Tz = (1 - z) /(1 + z)$. Recall that $T$ sends $\Pi^+$ to the unit disk $\D$, and therefore 
$f (z)\in \Pi^+$, for every $z \in  \Pi^+$ because $|\arg (\exp w)| < 1 <\pi/2$, for all $w \in \D$. \par

Now let $Q$ be the unit square $Q = (- 1, 1) \times (- 1, 1)$. For $0 < t\le 1/2$, let $Q_t$ be the square, centered in $1$, 
$Q_t = (1 - t, 1 + t) \times (- t, t)$, which is contained in $\Pi^+$, and define
\begin{displaymath}
\sigma (t) = \frac{1}{A (Q_t)} \int_{Q_t}|f (z)| \, dA(z) = \frac{1}{4 t^2} \int_{- t}^{t} \bigg[ \int_{1 - t}^{1 + t} |f (x+iy)| \, dx \bigg] \, dy\,.
\end{displaymath}
Using a change of variable we have:
\begin{displaymath}
\sigma (t) = \frac{1}{4} \iint_Q |f (1 + tx+ity)| \, dx \,dy \,.
\end{displaymath}

We are going to prove that there exists $t$ such that $\sigma (t) < 1 = |f(1)|$ and the average of $|f|$ in the cube $Q_t$ is smaller than $|f(1)|$.
Now observe that 
\begin{displaymath}
f (z) = \frac{1}{16} (z - 1)^4 +  O\, \bigl( (z - 1)^5 \bigr), \qquad z \to 1. 
\end{displaymath}
Consequently, there exists a constant $C > 0$ such that, for $z \in Q_{1/2}$,
\begin{displaymath}
\Re \bigl( f (z) \bigr) \le  \frac{1}{16} \, \Re\bigl( (z - 1)^4 \bigr) + C|z - 1|^5
\end{displaymath}
and then, there exists $C_1 > 0$, such that for every $x + i y \in Q$ and $t \in (0, 1/2)$, 
\begin{align*}
|f (1 + tx + i ty)|
& \le \exp\Bigl[ \frac{1}{16} \, \Re \bigl( t^4 (x + i y)^4 \bigr) + C_1 t^5 \Bigr] \\ 
& = \exp\Bigl[ \frac{t^4}{16} \, (x^4 + y^4 - 6 x^2 y^2) + C_1 t^5 \Bigr] \,.
\end{align*}
Integrating over $Q$, putting: 
\begin{displaymath}
\tau (s) = \frac{1}{4} \int_{- 1}^1 \! \int_{- 1}^1 \exp \bigl( (s/16) (x^4 + y^4 - 6 x^2 y^2) + C_1 s^{5/4} \bigr) \, dx \,dy\,,
\end{displaymath}
we get that $\sigma (t) \leq \tau (t^4)$, for $t \in (0, 1/2]$. We just need to prove that, for $s > 0$ close enough to $0$, we have $\tau (s) < 1$. But this is easy 
because $\tau (0) = 1$, and 
\begin{displaymath}
\tau ' (0) = \frac{1}{4} \int_{- 1}^1 \! \int_{- 1}^1 \frac{1}{16} \,  (x^4 + y^4 - 6 x^2 y^2) \, dx \, dy = 
\frac{1}{64} \Bigl( \frac{4}{5} + \frac{4}{5} - \frac{8}{3} \Bigr) =  - \frac{1}{60} < 0\,.
\end{displaymath}
\par
\bigskip\goodbreak

Return now to the proof of Theorem~\ref{theo clef}. \par

Consider, for every $n \geq 0$, the conditional expectation of the restriction to $\Omega$ of $|f|$ with respect to the algebra ${\cal Q}_n$ generated by 
the squares $Q_{(n, j, k)}$, $0 \leq j, k \leq 2^n - 1$ (note that ${\cal Q}_n \subseteq {\cal Q}_{n + 1}$): 
\begin{equation} 
(\esp_n |f|) (z)  = \sum_{j, k = 0}^{2^n - 1} \bigg( \frac{1}{A (Q_{(n, j, k)})} 
\int_{Q_{(n, j, k)}} |f|\, dA \bigg) \, \ind_{Q_{(n, j, k)}} (z) \,,
\end{equation} 
and the maximal function $M f$ is defined by:
\begin{equation} 
M f (z) = \sup_n \, (\esp_n |f|) (z)\,.
\end{equation}
One has
\begin{equation} 
M (f) (z)  = \sup_{z \in Q_{(n, j, k)}} \frac{1}{A (Q_{(n, j, k)})} \int_{Q_{(n, j, k)}} |f|\, dA \,.
\end{equation} 
Since $f$ is continuous on $\Omega$, one has $\lim_{n \to \infty} \esp_n |f| (z) = |f (z) |$ for every $z \in \Omega$, and it follows that:
\begin{equation} \label{inclusion} 
\{|f| > 1\} \subseteq \{M f > 1\} \,.
\end{equation} 
Now, the set $\{M f > 1\} \cap \Omega$ can be split into a disjoint union 
\begin{displaymath} 
\smash {\{M f > 1\} \cap \Omega = \bigsqcup_{n \geq 1} Z_n \,, } 
\end{displaymath} 
where 
\begin{displaymath} 
Z_n = \{z \in \Omega \, ;\ (\esp_n |f|) (z) > 1 \text{ and } (\esp_j |f|) (z) \leq 1 \text{ if } j < n \}  \,.
\end{displaymath} 
(note that, by Proposition~\ref{no beans},  $\esp_0 |f| \leq 1$ if $| f (1)|$ is small enough). \par 
\smallskip

Since $\esp_n |f|$ is constant on the sets $Q \in {\cal Q}_n$, each $Z_n$ can be in its turn decomposed, up to a set of measure $0$, into a disjoint union 
$E_n = \bigsqcup_{(j, k) \in J_n} Q _{(n, j, k)}$.\par
By definition, for $z \in Z_n$, one has $(\esp_n |f|) (z) \geq 1$ and hence, for $(j, k) \in J_n$,  
\begin{displaymath}
\qquad \frac {1} {A (Q_{(n, j, k)})} \int_{Q_{(n, j, k)}} |f| \, dA \geq 1 \qquad \text{for } z \in Q_{(n, j, k)} \,.
\end{displaymath}

But, on the other hand, $(\esp_{n - 1} |f|) (z) \leq 1$ for $z \in Z_n$, and we have, if $z \in Q_{(n, j, k)}$:
\begin{align*}
(\esp_n |f|) (z) 
& = \frac {1} {A (Q_{(n, j, k)})} \int_{Q_{(n, j, k)}} |f| \, dA 
\leq \frac {1} {A (Q_{(n, j, k)})} \int_{Q_{(n - 1, j', k')}} |f| \, dA  \\
& \leq 4\, \frac {1} {A (Q_{(n - 1, j', k')})} \int_{Q_{(n - 1, j', k')}} |f| \,dA  \leq 4\,,
\end{align*}
where $Q_{(n - 1, j', k')}$ is the square of rank $(n - 1)$ containing $Q_{(n, j, k)}$.\par
\smallskip

Finally, we can write $\{ M f > 1 \} \cap \Omega$ as a disjoint union, up to a set of measure $0$,  
\begin{equation}\label{union disjointe}
\{ M f > 1 \} \cap \Omega = \bigsqcup_{l \in L} Q_l \,, 
\end{equation}
where $L$ is a subset of all the indices $(n, j, k)$, for which:
\begin{equation}\label{encadrement moyennes}
1 \leq \frac {1} {A (Q_l)} \int_{Q_l} |f| \, dA \leq 4 \,.
\end{equation}

Equations \eqref{inclusion}, \eqref{union disjointe} and \eqref{encadrement moyennes} define the Calder{\'o}n-Zygmund decomposition of the function 
$f$. \par
\bigskip

We are now ready to end the proof of Theorem~\ref{theo clef}.\par
\medskip

For $\lambda \geq 1$, set $E_\lambda = \{ |f| > \lambda \}$; one has, by \eqref{union disjointe}, Proposition~\ref{auxiliary} and \eqref{moyenne Q}:
\begin{align*} 
\tau_\alpha (E_\lambda \cap \Omega) 
& = \tau_\alpha (E_\lambda \cap \{ Mf > 1\} \cap \Omega) 
= \sum_{l \in L} \tau_\alpha ( E_\lambda  \cap Q_l) \\
& \leq \frac{K_\alpha}{\lambda^{\alpha + 2}} \sum_{l \in L} \tau_\alpha (Q_l) \, |f (c_l)|^{\alpha + 2} \\
& \leq \frac{K_\alpha}{\lambda^{\alpha + 2}} \sum_{l \in L} \tau_\alpha (Q_l) \, \Big( \frac{16}{\pi} \Big)^{\alpha + 2} 
= \frac{C_\alpha}{\lambda^{\alpha + 2}} \sum_{l \in L} \tau_\alpha (Q_l)  \\
& = \frac{C_\alpha}{\lambda^{\alpha + 2}} \tau_\alpha (\{Mf > 1\} \cap \Omega) .
\end{align*} 
But, on the other hand, the sets $R_l$ of Lemma~\ref{precision} are disjoint, since $R_l \subseteq Q_l$ and we have 
$\vert f \vert > \delta_0 \, \vert f (c_l) \vert >  (4/\pi C) \, \delta_0 := \delta_1$ on $R_l$, in view of Lemma~\ref{precision} and  
Proposition~\ref{no beans}. Therefore:
\begin{align*} 
\tau_{\alpha} (\vert f \vert > \delta_1) 
& \geq \tau_{\alpha} \Big( \bigsqcup_{l} R_l \Big) 
= \sum_{l \in L} \tau_{\alpha} (R_l) 
\geq c \sum_{l \in L} \tau_{\alpha} (Q_l) 
= c \, \tau_{\alpha} \Big(\bigsqcup_{l \in L} Q_l \Big) \\
& = c \, \tau_{\alpha} (\{M f > 1\} \cap \Omega).
\end{align*} 
We get hence 
$\tau_\alpha (E_\lambda \cap \Omega) \leq \frac{C'_\alpha}{\lambda^{\alpha + 2}}\,  \tau_\alpha (\{|f| > \delta_1\})$ for $\lambda \geq 1$, 
with $C'_\alpha = C_\alpha/c$. Applying this to $f/\delta_1$ instead of $f$, we get:
\begin{displaymath} 
\tau_\alpha (\{|f| > \lambda\} \cap \Omega) \leq \frac{C{''}_\alpha}{\lambda^{\alpha + 2}}\,  \tau_\alpha (\{|f| > 1\})
\end{displaymath} 
for $\lambda > \lambda_1 := 1/ \delta_1$, and that finishes the proof of Theorem~\ref{theo clef}. 
\goodbreak


\section{An application to composition operators} 

In this section, we give an application of our main result to composition operators on weighted Bergman-Orlicz spaces.\par\medskip

Recall that an Orlicz function $\Psi \colon [0, \infty) \to \R_+$ is a non-decreasing convex function such that $\Psi (0) = 0$  
and $\Psi (x)/x \to \infty$ as $x$ goes to $\infty$. The \emph{weighted Bergman-Orlicz space} ${\mathfrak B}^\Psi_\alpha$ 
is the space of all analytic functions $f \colon \D \to \C$ such that 
\begin{displaymath} 
\int_\D \Psi (|f|/C)\, d{\cal A}_\alpha < +\infty 
\end{displaymath} 
for some constant $C > 0$. The norm of $f$ in ${\mathfrak B}^\Psi_\alpha$ is the infimum of the constants $C$ for which the above integral is 
$\leq 1$. With this norm, ${\mathfrak B}^\Psi_\alpha$ is a Banach space. \par
Now, every analytic self-map $\phi \colon \D \to \D$ defines a bounded linear operator 
$C_\phi \colon {\mathfrak B}^\Psi_\alpha \to {\mathfrak B}^\Psi_\alpha$ by $C_\phi (f) = f \circ \phi$, called the \emph{composition operator of 
symbol $\phi$}. This is a consequence of the classical Littlewood's subordination principle, using the facts that the measure ${\cal A}_\alpha$ is 
radial and the function $\Psi (|f|/C)$ is sub-harmonic for every analytic function $f \colon \D \to \C$. Such an operator may be seen as a 
Carleson embedding $J_\mu \colon {\mathfrak B}^\Psi_\alpha \to L^\Psi (\mu)$ for the pull-back measure $\mu = \phi ({\cal A}_\alpha)$. S. Charpentier 
(\cite{Charp}), following \cite{Bergman}, has characterized the compactness of such embeddings (actually in the more general setting of the unit ball 
${\mathbb B}_N$ of $\C^N$ instead of the unit disk $\D$ of $\C$):
\begin{theorem} [S. Charpentier]
For every finite positive measure $\mu$ on $\D$ and for every $\alpha > -1$, one has:\par
1) If ${\mathfrak B}^\Psi_\alpha$ is compactly contained in $L^\Psi (\mu)$, then 
\begin{equation} \label{cond necess}
\lim_{h \to 0} \frac{\Psi^{-1} (1/h^{\alpha + 2})}{\Psi^{-1} (1/\rho_\mu (h))} = 0\,.
\end{equation} 
\par
2) Conversely, if
\begin{equation} \label{cond suff}
\lim_{h \to 0} \frac{\Psi^{-1} (1/h^{\alpha + 2})}{\Psi^{-1} (1/ h^{\alpha + 2} K_\mu (h))} = 0\,,
\end{equation} 
then ${\mathfrak B}^\Psi_\alpha$ is compactly contained in $L^\Psi (\mu)$. 
\end{theorem} 

Here $\rho_\mu$ is the \emph{Carleson function} of $\mu$, defined as:
\begin{equation} \label{rho}
\rho_\mu (h) = \sup_{|\xi|=1} \mu [ W (\xi, h)] 
\end{equation} 
and 
\begin{equation} \label{K}
K_\mu (h) = \sup_{0 < t \leq h} \frac{\rho_\mu (t)}{t^{\alpha + 2}} \,\cdot
\end{equation} 

When $\mu = \phi ({\cal A}_\alpha)$ is the pull-back measure of ${\cal A}_\alpha$ by an analytic self-map $\phi \colon \D \to \D$, we denote 
them by $\rho_{\phi, \alpha + 2}$ and $K_{\phi, \alpha + 2}$ respectively. \par
\medskip

We gave in \cite{Bergman}, in the non-weighted case, examples showing that conditions \eqref{cond necess} and \eqref{cond suff} are not equivalent for 
general measures $\mu$. However, Theorem~\ref{theo principal} implies that 
$K_{\phi, \alpha + 2} (h) \lesssim \rho_{\phi, \alpha + 2} (h)/h^{\alpha + 2}$ and so conditions~\eqref{cond necess} and \eqref{cond suff} are 
equivalent in this case. Therefore, we get:

\begin{theorem} 
For every $\alpha > -1$, every Orlicz function $\Psi$, and every analytic self-map $\phi \colon \D \to \D$, the composition operator 
$C_\phi \colon {\mathfrak B}^\Psi_\alpha \to {\mathfrak B}^\Psi_\alpha$ is compact if and only if:
\begin{equation} 
\lim_{h \to 0} \frac{\Psi^{-1} (1/h^{\alpha + 2})}{\Psi^{-1} (1/\rho_{\phi, \alpha + 2} (h))} = 0\,.
\end{equation} 
\end{theorem} 
%



\end{document}